\documentclass[12pt]{article}
\input epsf

\usepackage{amsfonts}
\usepackage{amscd}
\usepackage{amsmath,amssymb,amsthm}
\usepackage{amstext}
\usepackage{amscd}
\usepackage{graphicx}
\usepackage[all]{xy}
\usepackage{pstricks,pst-node,pst-tree}

\theoremstyle{plain}
\newtheorem{lem}{Lemma}[section]
\newtheorem*{class}{Theorem}

\newtheorem{theo}[lem]{Theorem}

\theoremstyle{definition}

\newtheorem{definition}[lem]{Definition}

\newtheorem{rem}[lem]{Remark}

\renewcommand{\descriptionlabel}[1]%
       {\hspace{\labelsep}\textsf{#1}}

\newcommand{\bH}{\mathbb{bH}}
\newcommand{\R}{\mathbb{R}}

\newcommand{\noi} {\noindent}

\DeclareMathOperator{\Int}{Int}

\begin{document}
\title{Topological surfaces as gridded surfaces in geometrical spaces.}
\author{ Juan Pablo D\'iaz\thanks{This work was supported by FORDECYT 265667 (M\'exico)},
 Gabriela Hinojosa, Alberto Verjosvky, \thanks{This work was partially supported by
CONACyT (M\'exico), CB-2009-129280  and PAPIIT (Universidad
Nacional Aut\'onoma de M\'exico) \#IN106817.}}

\date{November 27, 2017}
\maketitle

\begin{abstract} 
\noi In this paper we study topological surfaces as gridded surfaces in the 2-dimensional scaffolding of cubic honeycombs in Euclidean and hyperbolic spaces. 
\end{abstract}

\noi {\bf Keywords:} Cubulated surfaces, gridded surfaces, euclidean and hyperbolic honeycombs. 

\noi  {\bf AMS subject classification:} Primary 57Q15, Secondary 57Q25, 57Q05.

\section{Introduction}

The category of cubic complexes and cubic maps is similar to the simplicial category. The only
difference consists in considering cubes of different dimensions instead of simplexes. In this context,
a \emph{cubulation} of a manifold  is a cubical complex which is PL homeomorphic to the
manifold (see \cite{dolbilin}, \cite{funar}, \cite{matveev}). In this paper we study the realizations of  cubulations of manifolds embedded in skeletons (or scaffoldings) of the canonical cubical honeycombs of an euclidean or hyperbolic space.\\ 

\noindent In \cite{BHV} it  was shown the following theorem:

\begin{theo} Let  $M,\,\,N\subset\Bbb R^{n+2}$, $N\subset{M}$,
be closed and smooth submanifolds of  $\Bbb R^{n+2}$ such that $\mbox{dimension}\,(M)=n+1$ and $\mbox{dimension}\,(N)=n$.
Suppose that $N$ has a trivial normal bundle in $M$ ({\it i.e.,} $N$ is a two-sided hypersurface of $M$). Then there exists an ambient isotopy of 
$\Bbb R^{n+2}$ which takes $M$ into the $(n+1)$-skeleton of the canonical
cubulation $\cal C$ of $\Bbb R^{n+2}$ and $N$ into the $n$-skeleton of $\cal C$. In particular, $N$ can be deformed by
an ambient isotopy into a cubical manifold contained in the canonical scaffolding of  $\Bbb R^{n+2}$.
\end{theo}

\noindent In particular, the previous theorem establishes that smooth knotted surfaces in $\Bbb R^{4}$ are isotopic to cubulated 2-knots in 
the 4-dimensional cubic honeycomb. Orientable smooth closed surfaces are cubulated manifolds.

\begin{definition}
Let $S$ be a topological $2$-manifold embedded in either $\mathbb{R}^{n}$ or $\mathbb{H}^{n}$. We say that $S$ is a \emph{geometrically gridded surface} 
if it is contained in the
$2$-skeleton of either the canonical cubulation  $\{4,3^{n-2},4\}$ of  $\mathbb R^{n}$,  the hyperbolic cubic honeycomb $\{4,3,5\}$ of the hyperbolic 
space $\mathbb{H}^{3}$
or the hyperbolic cubic honeycomb  $\{4,3,3,5\}$ of  the hyperbolic space $\mathbb{H}^{4}$. Here \emph{geometrically} means
that the surface is gridded by isometric pieces of regular euclidean or hyperbolic squares and it is placed in the corresponding skeleton in the 
scaffolding $\mathcal{S}^2$ of a regular cubic euclidean or hyperbolic honeycomb $\cal C$. If it is understood the type of squares we simply say that 
$S$ is a \emph{gridded surface}.
\end{definition}

\noi A \emph{gridded surface} $S$ on $\mathcal{S}^2$ of $\cal C$ in $\mathbb R^{n}$ is a 
piecewise linear surface such that each linear piece is a unit square with its vertices in the $\mathbb{Z}^n$-lattice of $\mathbb R^{n}$. 
However in the case of hyperbolic spaces $\bH^n$
the description of the vertices is more complicated since the vertices belong to an orbit of a non-abelian discrete group acting by isometries on hyperbolic space 
$\bH^n$ as we will see in the next sections. 
We remark that our gridded surfaces are in a natural way length spaces \cite{BB}.\\

\noi Hilbert (1901, \cite{hilbert}) proved that there is no regular smooth isometric immersion $X :\mathbb{H}^2\rightarrow\mathbb{R}^3$.  
Efimov (1961, \cite{efimov}) generalizad this nonexistence theorem of Hilbert to the case of complete surfaces of nonpositive curvature; more precisely, he showed that there is no ${\cal{C}}^2$-isometric immersion of a complete, two dimensional, Riemannian manifold 
$M\subset\mathbb{R}^3$ whose curvature satisfies $K\leq c<0$. However, J. Nash (1956, \cite{struik}) proved that any ${\cal{C}}^k$-manifold $(M^n,g)$ 
can be ${\cal{C}}^k$-isometrically immersed into $\mathbb{R}^n$ where $q\geq \frac{3}{2}n(n+1)(n+9)$ and $k\geq 3$. This implies that 
the hyperbolic space can be embedded into a high dimensional Euclidean space; for instance, we can find a  ${\cal{C}}^3$-isometric embedding 
of $\mathbb{H}^2$ into $\mathbb{R}^{99}$ (see \cite{gromov}, \cite{struik}). \\

\noi In this paper we study topological surfaces as geometrically gridded surfaces in the scaffolding of cubic honeycombs in Euclidean and hyperbolic spaces.
We prove that connected orientable surfaces can be gridded in $\mathbb{R}^{3}$ or  $\mathbb{H}^{3}$ and all the surfaces, orientable or not, can be gridded in $\mathbb{R}^{4}$ or  $\mathbb{H}^{4}$. 

\section{Preliminaries}
\noi This section consists of two topics. The first studies the type of ``scaffolding'' (\emph{i.e.}. the 2-skeleton of an 
Euclidean or hyperbolic honeycomb) in which we can embed a surface in order to make it a gridded surface. For this purpose we start studying the 
regular cubic honeycombs of dimensions 3 and 4 in the Euclidean and hyperbolic cases. 

\noi In the second topic we will revise the theorem of topological classification of non-compact surfaces, in particular non-compact surfaces with
Cantor sets of ends of planar and nonplanar type, and also with Cantor sets of non-orientable ends.    

\subsection{Regular cubic honeycombs}

\noindent  We are interested in geometrically regular cubic  honeycombs which are geometric spaces filled with hypercubes which are euclidean or hyperbolic hypercubes. We denote these honeycombs by their Schl\"affli symbols. For a cube the symbol is $\{4,3\}$. This means that the faces of the regular cube are squares with Schl\"affli symbol  $\{4\}$ and that there are 3 squares around each vertex. These cubic honeycombs have Schl\"affli symbols which describe their geometry and start by $\{4,3,...\}.$

\subsubsection{Euclidean cubic honeycombs $\{4,3^{n-2},4\}$}

\noindent  The {\it canonical cubulation} $\mathcal{C}^n$ of $\mathbb{R}^{n}$ is its decomposition into $n$-dimensional cubes which are the images of the unit $n$-cube
$I^n=[0,1]^n$ by translations by vectors with integer coefficients. Then all vertices of $\mathcal{C}^n$ have integers in their coordinates.\\

\noindent Any {\it cubulation} of $\mathbb{R}^{n}$ is obtained by applying a conformal transformation to the canonical cubulation. 
Remember that a {\it conformal transformation} is of the form  
$x\mapsto{\lambda{A(x)}+a},$ where $\lambda\neq0,\,\,a\in \mathbb{R}^{n},\,\,\,A\in{SO(n)}$.\\

\noi Any $n-$cubulation has the same combinatorial structure as honeycomb. The regular hypercubic honeycomb whose Schl\"affli symbol is $\{4,3^{n-2},4\}$ is a {\it cubulation} of $\mathbb{R}^{n}$ which is its 
decomposition into a collection $\mathcal{C}^n$ of right-angled $n$-dimensional
hypercubes $\{4,3^{n-2}\}$ called the $\textit{cells}$ such that any two are
either disjoint or meet in one common $k-$face of some dimension $k$. This provides  $\mathbb{R}^{n}$ with the structure of a cubic
complex whose category is similar to the simplicial category PL. \\

\noi The combinatorial structure of the regular honeycomb $\{4,3,4\}$ is as follows: there are 6 edges, 12 squares and 8 cubes which are incident 
for  each vertex and there are 4 squares and 4 cubes which are incident 
for  each edge. \\

\noi The combinatorial structure of the regular honeycomb $\{4,3,3,4\}$ is as follows: there are 8 edges, 24 squares, 32 cubes and 16 hypercubes which are incident 
for  each vertex; there are 6 squares, 32 cubes and 16 hypercubes which are incident 
for  each edge and there are 4 cubes and 4 hypercubes which are incident 
for  each square.\\  

\begin{figure}[h]  
\begin{center}
\includegraphics[height=5cm]{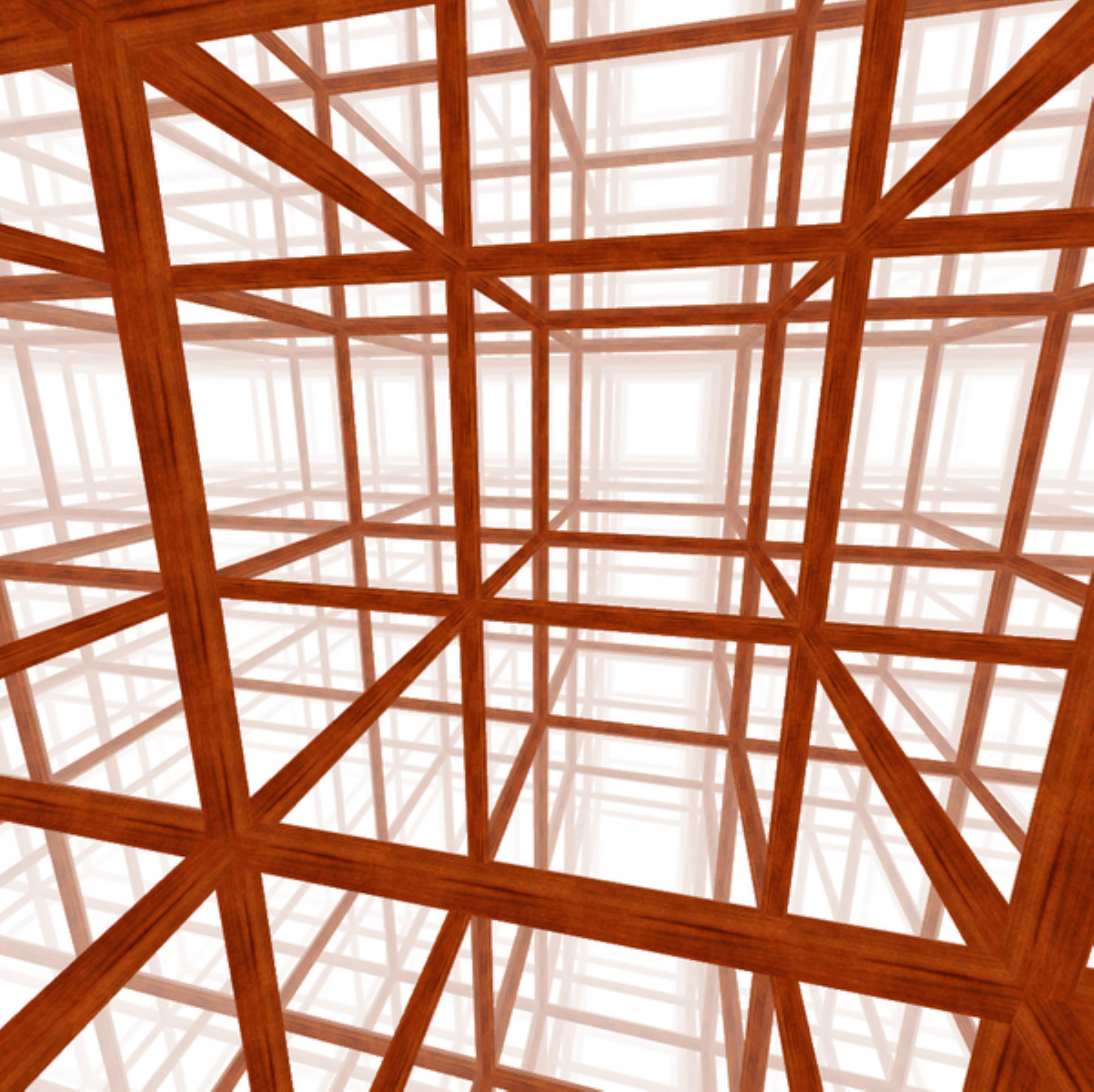}
\end{center}
\caption{\sl The 3-dimensional cubic kaleidoscopic honeycomb $\{4,3,4\}$.     This Figure is courtesy of Roice Nelson \cite{RN}.} 
\label{F1}
\end{figure}

\begin{definition} The $k$\emph{-skeleton} of the canonical cubulation $\mathcal{C}^n$ of $\mathbb{R}^{n}$, denoted by $\mathcal{C}^k$, 
consists of the union of the $k$-skeletons of the hypercubes in  $\mathcal{C}^n$,
{\it i.e.,} the union of all cubes of dimension $k$ contained in  the $n$-cubes in  $\mathcal{C}^n$.
We will call the 2-skeleton $\mathcal{C}^2$ of $\mathcal{C}^n$ the {\it canonical scaffolding} of $\mathbb{R}^{n}$.
\end{definition}

\subsubsection{Hyperbolic cubic honeycombs $\{4,3,5\}$ and $\{4,3,3,5\}$}

\noi A gridded surface is a surface made of congruent squares contained in the corresponding scaffoldings which have disjoint interiors and are glued only in their edges, two squares are either disjoint or share one edge or a vertex. Around each vertex there is a circular circuit of squares (the \emph{squared star} of the vertex). The condition of regularity of platonic solids in the 3-dimensional case implies that the number of squares around each vertex is a constant $k>2$. If $k=3$ then the regular gridded surface is the cube $\{4,3\}$ which is homeomorphic to a sphere $\mathbb{S}^2$. If $k=4$ then the regular surface is the gridded plane $\{4,4\}$ which is isometric to the Euclidean plane $\mathbb{R}^2$. If $k>4$ then the regular surfaces are the gridded hyperbolic planes which are denoted by $\{4,5\}, \{4,6\},...$. These gridded hyperbolic planes are length spaces \cite{BB} and are isometric, as length metric spaces to the hyperbolic plane $\mathbb{H}^2$. \\

\noi Analogously, we consider regular gridded 3-manifolds made of congruent cubes in scaffoldings which are disjoint and glued only in their boundary squares, two cubes are either disjoint or meet at a vertex an edge or a square.  There is a circular circuit of cubes around each edge and one spherical circuit of cubes around each vertex as PL 3-manifolds. For each edge there is the same number of cubes $k>2$. If $k=3$ then the regular gridded 3-manifold is the hypercube $\{4,3,3\}$ which is homeomorphic to the 3-sphere $\mathbb{S}^3$. If $k=4$ then the regular gridded 3-manifold is the cubic space $\{4,3,4\}$ which is isometric to the Euclidean space $\mathbb{R}^3$. If $k=5$ then the regular gridded 3-manifold is the cubic hyperbolic 3-space $\{4,3,5\}$ which is isometric as a length space to the hyperbolic 3-space $\mathbb{H}^3$ \cite{BB}.  \\

\begin{figure}[h]  
\begin{center}
\includegraphics[height=5cm]{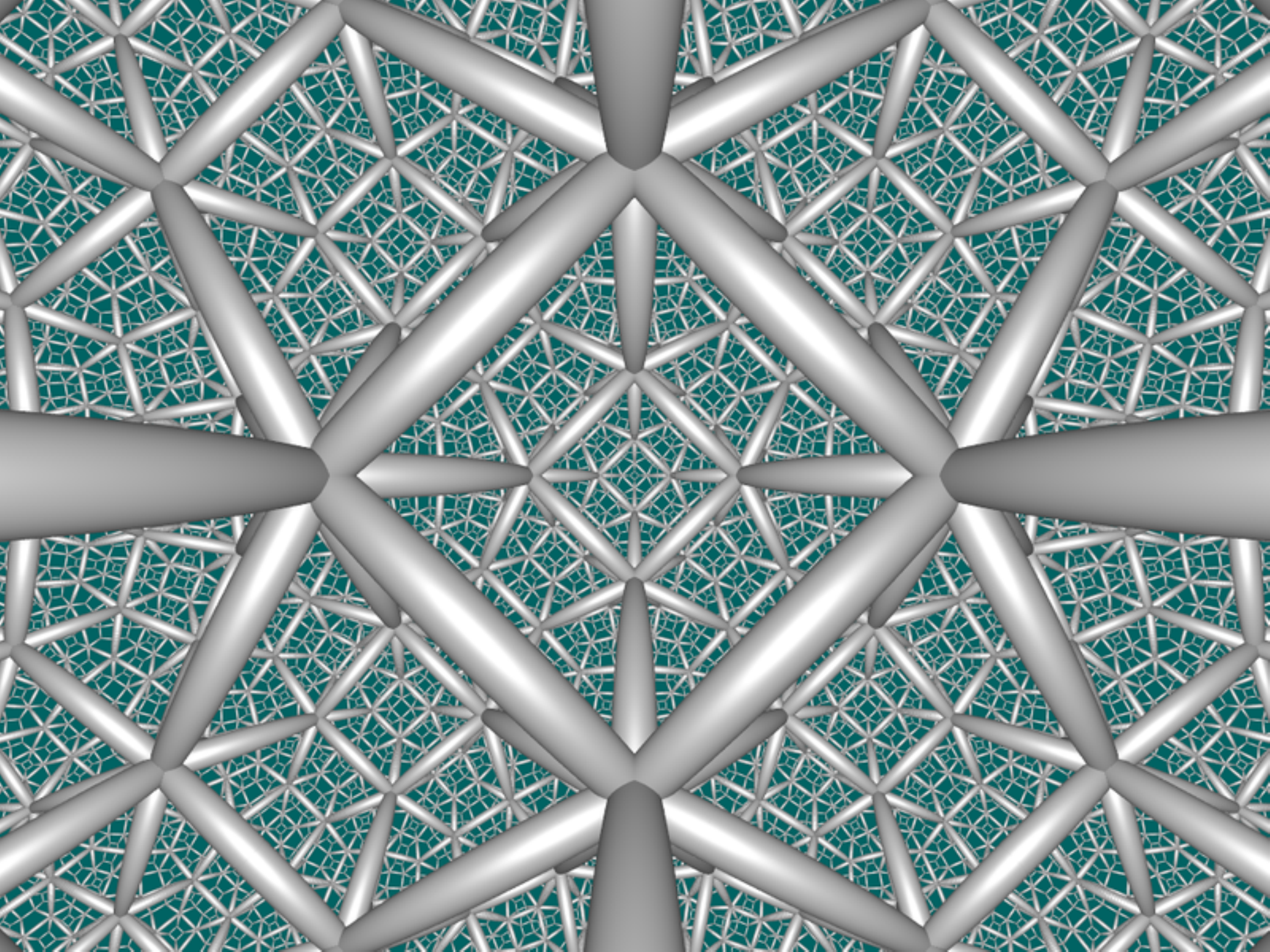}
\end{center}
\caption{\sl The 3-dimensional cubic hyperbolic honeycomb $\{4,3,5\}$. This Figure is courtesy of Roice Nelson \cite{RN}.} 
\label{F10}
\end{figure}

\noi Finally we construct regular gridded 4-manifolds made of congruent hypercubes contained in the corresponding scaffoldings which are disjoint and glued only in their boundary cubes, two hypercubes are either disjoint, or meet along a boundary cube of some dimension.  
There is one circular circuit of hypercubes around each square, one 2-spherical circuit of hypercubes around each edge and one 3-spherical circuit of hypercubes around each edge vertex as PL 4-manifolds. \\

\noi For each square which is a ridge, there is the same number  of hypercubes $k>3$. If $k=3$ then the regular gridded 4-manifold is the hypercube $\{4,3,3,3\}$ which is 
homeomorphic to the 4-sphere $\mathbb{S}^4$. If $k=4$ then the regular 4-manifold is the gridded cubic 4-space $\{4,3,3,4\}$ which is isometric to the Euclidean 4-space $\mathbb{R}^4$. 
If $k=5$ then the regular 4-manifold is the gridded hyperbolic 5-space $\{4,3,3,5\}$ which is isometric to the hyperbolic 4-space $\mathbb{H}^4$. \\

\noi The combinatorial structure of the regular hyperbolic honeycomb $\{4,3,5\}$ ($\cal C$) is as follows: there are 12 edges, 30 squares and 20 cubes meeting at every vertex and there are 5 squares and 5 cubes meeting at every edge. \\

\noi The combinatorial structure of the regular honeycomb $\{4,3,3,5\}$ ($\cal C$) is as follows: there are 120 edges, 720 squares, 1200 cubes and 600 hypercubes meeting at every vertex; 
there are 6 squares, 32 cubes and 16 hypercubes meeting at every edge and there are 5 cubes and 5 hypercubes meeting at every square.\\  

\noindent As before, we define the $k$\emph{-skeleton} of the hyperbolic honeycomb of $\mathbb{H}^{n}$ ($n=3,4$), denoted by $\mathcal{C}^k$, 
consists of the union of the $k$-skeletons of the hypercubes in  $\cal C$,
{\it i.e.,} the union of all cubes of dimension $k$ contained in  the $n$-cubes in  $\cal C$.
We will call the 2-skeleton $\mathcal{C}^2$ of $\cal C$ the {\it canonical scaffolding} of $\mathbb{H}^{n}$ ($n=3,4$).

\subsection{Taxonomy of topological surfaces}
Next we will consider all topological surfaces. First we will remind some definitions.\\

\noindent Let $g(S)$ denote the genus of a compact surface $S$. The 
\emph{genus} of the noncompact surface $S$ is by definition
$$g(S)=\max\{g(A) \,:\, A \,\,\mbox{is a compact subsurface of}\,\, S\},$$ if this maximum exists, and infinite otherwise ( $g(S)=\infty$). \\

\noindent If $g(S)=0$ we say that the surface $S$ is \emph{planar}. In this case the surface is homeomorphic to an open subset of the plane.

\begin{definition}
There are four orientability classes of noncompact surfaces:
\begin{enumerate}
\item  If every compact subsurface of a surface $S$ is orientable, then $S$ is \emph{orientable}.
\item  If there is no bounded subset $A$ of $S$ such that $S-A$ is orientable, then $S$ is \emph{infinitely
nonorientable}.
\item If $S$ does not belong to (l) or (2) and every sufficiently large subsurface of $S$ contains an even number of cross caps, then $S$ is \emph{even non orientable}.
\item If $S$ does not belong to (1) or (2) and every sufficiently large subsurface of $S$ contains an odd number of cross caps, then $S$ is \emph{odd non orientable}.
\end{enumerate}
\end{definition}

\begin{definition}
Let $S$ be a surface. An \emph{end} of $S$ is a function $e$ which asigns to each compact subset $K$ of $S$ an unbounded, connected component
$e(K)$ of $\mbox{Cl}(S- K)$ in such a way that if $K$ and $L$ are compact subsets of $S$ and $K\subset L$ then $e(L)\subset e(K)$.
$E_S$ denotes the set of all ends of $S$.
\end{definition}

\noindent If $S$ is a noncompact surface, there is a compact subsurface $A$ of $S$ such that each component of $\mbox{Cl}(S - A)$ is either orientable or 
infinitely nonorientable, and either planar or of infinite genus.  In all that follows $A$ will denote such a subsurface.

\begin{definition}
Let $e\in E_S$. We say that $e$ is nonorientable (or orientable) if $e(A)$ is infinitely nonorientable (or orientable). And $e$ is planar (or of infinite genus) 
if $e(A)$ is planar (or of infinite genus).
\end{definition}

\noindent For any surface $S$ consider the nested triple $(E_S,G_S, O_S)$, where $E_S$ is the set of ends of $S$, $G_S$ is the subset of $E_S$
 consisting of the ends which are not planar, and $O_S$ is the subset of $G_S$ of the orientable nonplanar ends. \\

\noindent The next Theorem is  by Ker\'ekj\'arto and Richards, the reader can find it in \cite{richards}.

\begin{class}\emph{(Classification theorem for noncompact surfaces)}. 
Let $X$ and $Y$ be two noncompact surfaces of the same genus and orientability class.
Then $X$ and $Y$ are homeomorphic if and only if the triads $(E_X,G_X,O_X)$ and $(E_Y,G_Y, O_Y)$ are topologically equivalent.
\end{class}

\begin{rem}
The condition that $X$ and $Y$ are of the same genus and orientability class guaranties that their bounded parts are homeomorphic. 
The condition on the subsets of ends makes sure that their asymptotic behavior is the same.
\end{rem}

\noi Our goal is to show that any surface is homeomorphic to a gridded surface. In order to prove it, 
we will use the classification theorem for noncompact surfaces and the theorems of Richards  which appear in \cite{richards}.

\begin{theo} The set of ends of a surface is totally disconnected,
separable, and compact. Any compact, separable, totally disconnected space is
homeomorphic to a subset of the Cantor set.
\end{theo}

\begin{theo} Let $(X, Y, Z)$ be any triple of compact, separable, totally disconnected
spaces with $Z\subset Y\subset X$. Then there is a surface S whose ideal boundary
$(E_S, G_S, O_S)$ is topologically
equivalent to the triple $(X, Y, Z)$.
\end{theo}

\begin{theo}\label{rich} Every surface is homeomorphic to a surface formed from a
2-sphere $\mathbb{S}^2$ by first removing a closed totally disconnected set $X$ from $\mathbb{S}^2$, then
removing the interiors of a finite or infinite sequence $D_1, D_2, \dots$ of non-overlapping
closed discs in $\mathbb{S}^2 - X$, and finally suitably identifying the boundaries of these
discs in pairs, (It may be necessary to identify the boundary of one disc with
itself to produce an odd "cross cap.'') The sequence $D_1,D_2, \dots$ ''approaches
$X$" in the sense that, for any open set $U\subset \mathbb{S}^2$ containing $X$, all but a finite
number of the $D_i$ are contained in $U$.
\end{theo}

\noi Let $X\subset \mathbb{S}^2$ a Cantor set, the orientable noncompact surface  $\mathbb{S}^2 - X$ is called \emph{the tree of life}. 
The \emph{tree of life} can be constructed by an infinite set of pair of pants glued along their boundaries. A pair of pants is
by definition the surface with boundary $\mathbb{S}^2 - \{D_1,D_2,D_3\}$ 
where each $D_i$ is an open disc. Consider the surfaces with boundary $\Sigma_1$  and $\Sigma_2$ which are the connected sum of a pair 
of pants with a torus $\Sigma_1=\mathbb{T}^2 - \{D_1,D_2,D_3\}$ and with a projective plane $\Sigma_2=\mathbb{P}^2 - \{D_1,D_2,D_3\}$. \\

\noi \emph{A pruned tree of life} is obtained from the tree of life removing a set of branches (neighborhoods of ends) cut along 
boundaries of the corresponding pair of pants where these boundaries are identified to points (i.e, we cap some holes of the pair of pants). \\

\noi For the purpose of this paper we will rephrase the Theorem \ref{rich}, as follows.

\begin{theo}\label{clas} Every surface is homeomorphic to a surface obtained from a pruned tree of life such that 
a finite or infinite number of pair of pants  are interchanged by either the connected sum of a pair of pants with a torus 
{\it i.e.,} $\Sigma_1$ or the connected sum of pair of pants with a projective plane {\it i.e.,}  $\Sigma_2$. 
\end{theo}

\begin{figure}[h]  
\begin{center}
\includegraphics[height=2cm]{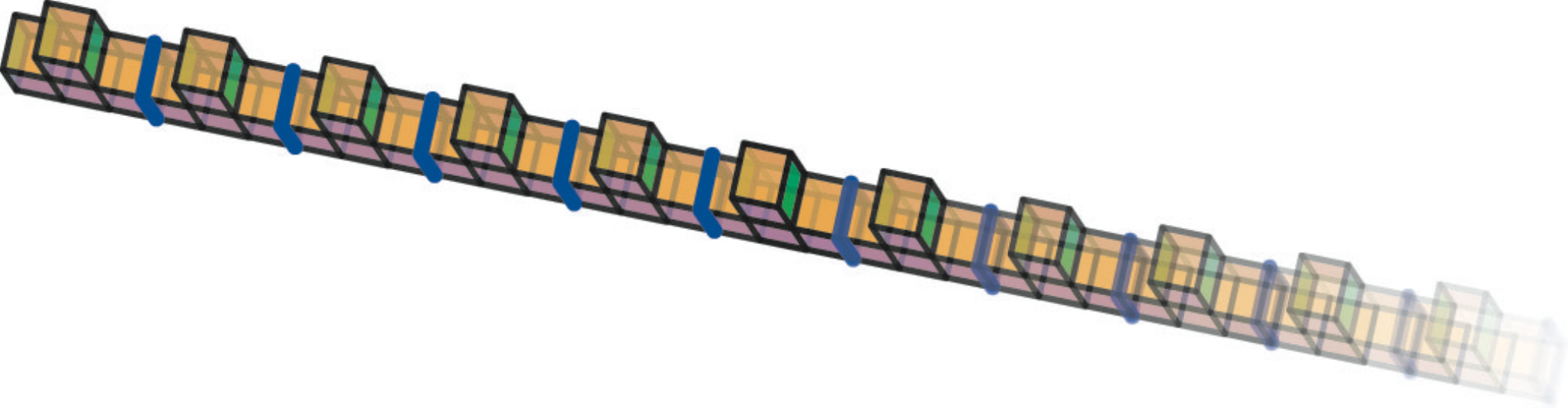}
\end{center}
\caption{\sl The gridded plane is obtained from the tree of life made of pair of pants after pruning and filling the remaining holes with squares (the left square and the squares at the top).} 
\label{cyl}
\end{figure}

\noi For instance, consider the \emph{plane}; \emph{i.e.} the noncompact surface  homeomorhic to 
$\mathbb{R}^2$. Then, it can be
obtained  from the tree of life if one removes all of its branches except one but and at each boundary component we glue a squared disk.  The gridded surface obtained 
in this way is depicted in Figure \ref{cyl}. 
In this way, the plane is composed by an infinite number of pair of pants.

\section{Surfaces in $\{4,3^{n-2},4\}$ in $\mathbb{R}^n$}

\subsection{Compact surfaces in $\{4,3,3,4\}$ in $\mathbb{R}^4$}

As we mentioned before, we will prove that any connected surface with a finite number of ends is homeomorphic to a gridded surface. We will start 
considering closed surfaces. \\

\begin{figure}[h]  
\begin{center}
\includegraphics[height=3cm]{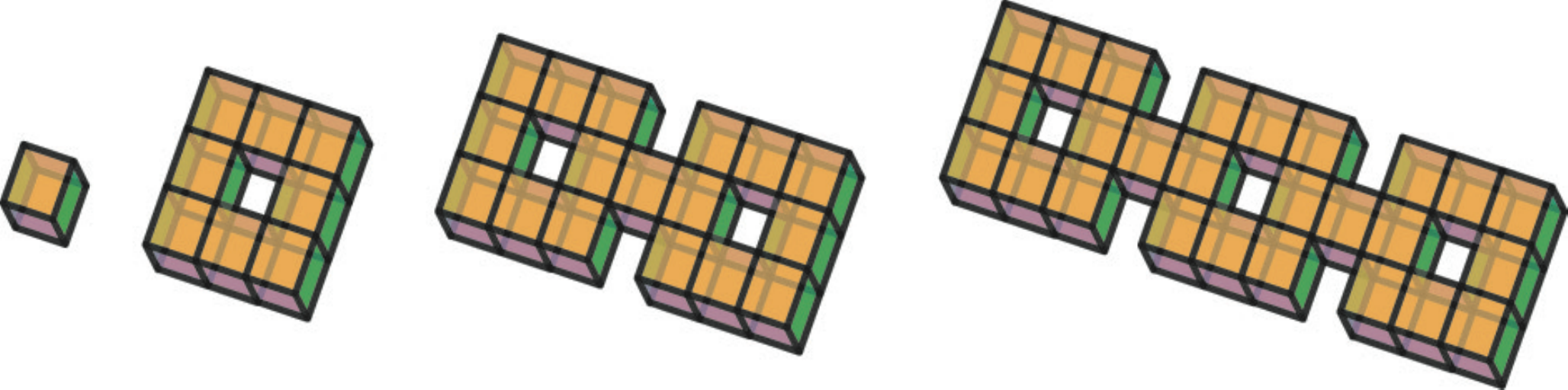}
\end{center}
\caption{\sl Closed gridded surfaces in $\mathbb{R}^3$.} 
\label{F2}
\end{figure}

\noindent Notice that the 2-sphere $\mathbb{S}^2$ and the 2-torus $\mathbb{T}^2$ are homeomorphic to gridded surfaces contained in 
the scaffolding $\mathcal{C}^2$ of the canonical cubulation $\cal C$  of $\mathbb{R}^{3}$ in the obvious way. In fact,
Consider the unit 3-cube $I^3=[0,1]^3$, clearly its boundary $\partial I^3$ is contained in $\mathcal{C}^2$ and
is homeomorphic to the 2-sphere $\mathbb{S}^2$ (see Figure \ref{F2}).\\

\noindent The 2-torus $\mathbb{T}^2$ is the gridded surface $S=\cup_{i=0}^9 F_i$ (see the second image of the Figure \ref{F2}), where  $F_i$, $i=0,\ldots,9$ are the following 10 sets which are unions of squares (squared sets): 
$$
F_0=\{(x,y,0)\in\mathbb{R}^3\,|\,0\leq x,\,y\leq 3\}\setminus \{(x,y,0)\in\mathbb{R}^3\,|\,1< x,\,y < 2\},
$$
$$
F_1=\{(x,y,1)\in\mathbb{R}^3\,|\,0\leq x,\,y\leq 3\}\setminus \{(x,y,0)\in\mathbb{R}^3\,|\,1< x,\,y < 2\},
$$
$$
\hspace{-3.8cm}F_2=\{(0,y,z)\in\mathbb{R}^3\,|\,0\leq y\leq 3,\,\,0\leq z\leq 1\},
$$
$$
\hspace{-3.8cm}F_3=\{(1,y,z)\in\mathbb{R}^3\,|\,0\leq y\leq 3,\,\,0\leq z\leq 1\},
$$
$$
\hspace{-3.8cm}F_4=\{(x,0,z)\in\mathbb{R}^3\,|\,0\leq x\leq 3,\,\,0\leq z\leq 1\},
$$
$$
\hspace{-3.8cm}F_5=\{(x,1,z)\in\mathbb{R}^3\,|\,0\leq x\leq 3,\,\,0\leq z\leq 1\},
$$
$$
\hspace{-3.8cm}F_6=\{(x,1,z)\in\mathbb{R}^3\,|\,1\leq x\leq 2,\,\,0\leq z\leq 1\},
$$
$$
\hspace{-3.8cm}F_7=\{(x,2,z)\in\mathbb{R}^3\,|\,1\leq x\leq 2,\,\,0\leq z\leq 1\},
$$
$$
\hspace{-3.8cm}F_8=\{(1,y,z)\in\mathbb{R}^3\,|\,1\leq y\leq 2,\,\,0\leq z\leq 1\},
$$
$$
\hspace{-3.8cm}F_9=\{(2,y,z)\in\mathbb{R}^3\,|\,1\leq y\leq 2,\,\,0\leq z\leq 1\}.
$$ 

\begin{rem}
 Notice that each square $F$ of the canonical cubulation $\cal C$ of $\mathbb{R}^{4}$ (or $\mathbb{R}^{3}$) is determined by its barycenter $B_F$. In fact, consider the unitary canonical vectors on $\mathbb{R}^4$: 
$e_{\pm 1}=(\pm 1,0,0,0)$, $e_{\pm 2}=(0,\pm 1,0,0)$, 
$e_{\pm 3}=(0,0,\pm 1,0)$ and $e_{\pm 4}=(0,0,0,\pm 1)$. Then
$$
F=\{ae_u+be_v\,:\,0\leq a,\,b\leq 1\}+w
$$ 
where $e_u$ and $e_v$ ($u,v\in\{\pm 1,\,\pm 2,\,\pm 3,\,\pm 4\}$, $|u|\neq |v|$),  denote the corresponding unitary canonical vectors and $w$ is a vector with
integers in its coordinates, in fact $w$ is a translation vector. Thus $B_F=\frac{1}{2}e_u+\frac{1}{2}e_v+w$.
\end{rem}

\noindent We identify the squares of the torus with their barycenters, then:
 
\noindent $F_0=(\frac{1}{2},\frac{1}{2},0), (\frac{3}{2},\frac{1}{2},0), (\frac{5}{2},\frac{1}{2},0), (\frac{1}{2},\frac{3}{2},0), (\frac{5}{2},\frac{3}{2},0), (\frac{1}{2},\frac{5}{2},0), (\frac{3}{2},\frac{5}{2},0), (\frac{5}{2},\frac{5}{2},0).$

\noindent $F_1=(\frac{1}{2},\frac{1}{2},1), (\frac{3}{2},\frac{1}{2},1), (\frac{5}{2},\frac{1}{2},1), (\frac{1}{2},\frac{3}{2},1), (\frac{5}{2},\frac{3}{2},1), (\frac{1}{2},\frac{5}{2},1), (\frac{3}{2},\frac{5}{2},1), (\frac{5}{2},\frac{5}{2},1).$

\noindent $F_2=(0,\frac{1}{2},\frac{1}{2}), (0,\frac{3}{2},\frac{1}{2}), (0,\frac{5}{2},\frac{1}{2}).$

\noindent $F_3=(1,\frac{3}{2},\frac{1}{2}).$

\noindent $F_4=(2,\frac{3}{2},\frac{1}{2}).$

\noindent $F_5=(3,\frac{1}{2},\frac{1}{2}), (3,\frac{3}{2},\frac{1}{2}), (3,\frac{5}{2},\frac{1}{2}).$

\noindent $F_6=(\frac{1}{2},0,\frac{1}{2}), (\frac{3}{2},0,\frac{1}{2}), (\frac{5}{2},0,\frac{1}{2}).$

\noindent $F_7=(\frac{3}{2},1,\frac{1}{2}).$

\noindent $F_8=(\frac{3}{2},2,\frac{1}{2}).$

\noindent $F_9=(\frac{1}{2},3,\frac{1}{2}), (\frac{3}{2},3,\frac{1}{2}), (\frac{5}{2},3,\frac{1}{2}).$\\

\noi The Klein bottle and the real projective plane can not be homeomorphic to gridded surfaces contained in the 2-skeleton of $\mathbb{R}^{3}$ (see \cite{Hans}) but they are homeomorphic to gridded surfaces contained in the 2-skeleton of $\mathbb{R}^{4}$.\\
\begin{lem}
The projective plane $(\mathbb{P}^2)$ is a gridded surface in $\{4,3,3,4\}$ in $\mathbb{R}^4$. 
\end{lem}

\begin{figure}[h]  
\begin{center}
\includegraphics[height=5cm]{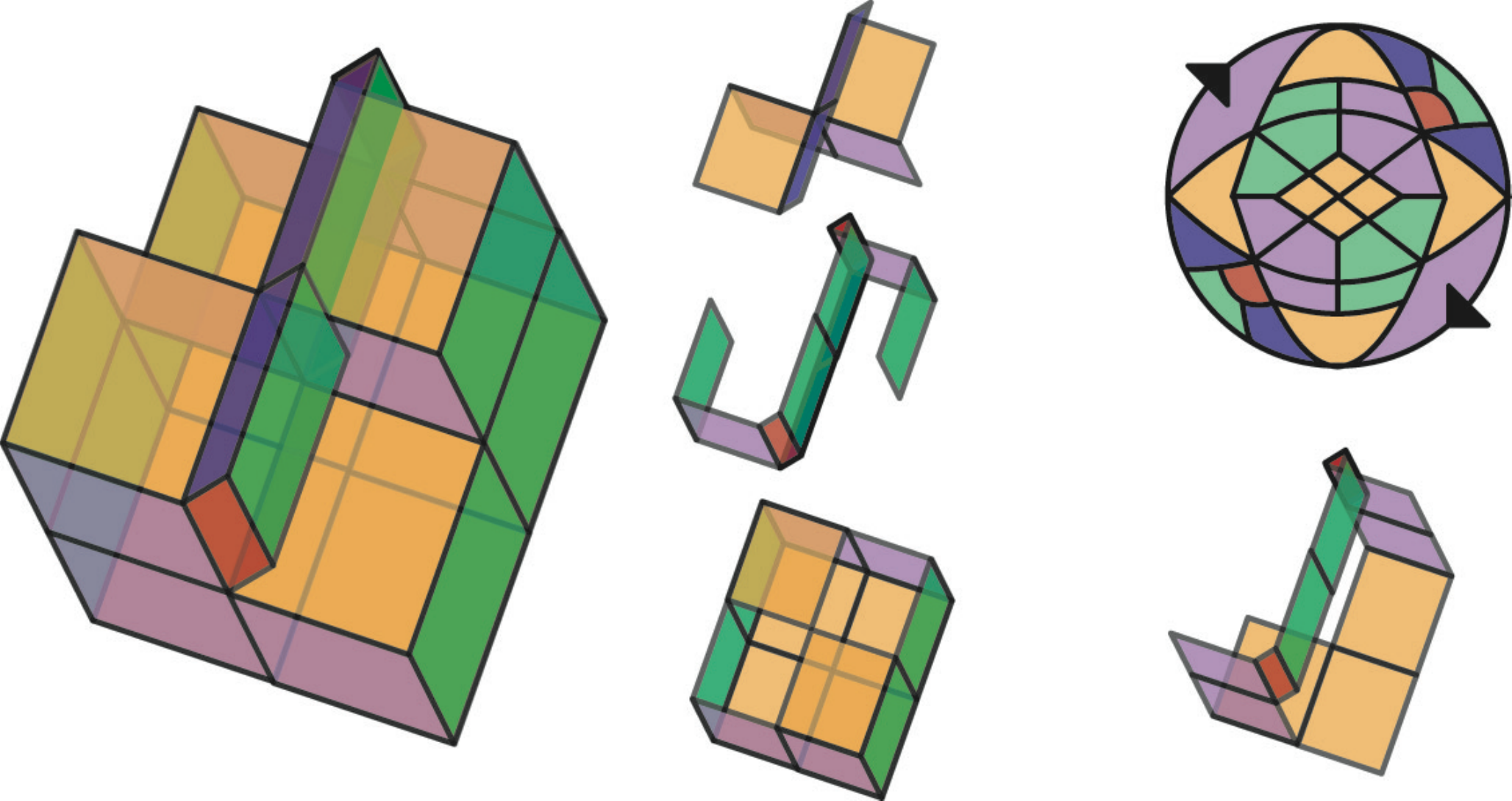}
\end{center}
\caption{\sl A gridded projective plane in $\mathbb{R}^4$. The first figure at the left is the projection of the 
crosscap into $\R^3$.} 
\label{PP}
\end{figure}

\noindent \emph{Proof}. We construct a gridded version of the crosscap in $\{4,3,3,4\}$ in $\mathbb{R}^4$, see Figure \ref{PP}. In the left we show the projection of the crosscap. 
In the middle we divide the crosscap in three parts: on the bottom, there is the base which is a cubic box minus two squares, in the middle there is a band and on the top, there is a disk which is the neighborhood of one vertex.    

\noi In the top right part of the Figure \ref{PP}, we can find the description of the combinatorial square complex of the crosscap as a  disk consisting on 30 squares, such that points in the circle boundary are identified by the antipodal map. 
On the bottom, there is a M\"obius band contained in this crosscap. 

\noindent The gridded crosscap is formed by 30 squares in planes parallel to five of the six coordinate planes in $\mathbb{R}^4$ and whose barycenters are:\\

\noindent $XY$: $(\frac{1}{2},\frac{1}{2},0,0),$ $(\frac{3}{2},\frac{1}{2},0,0),$ $(\frac{3}{2},\frac{3}{2},0,0),$ $(\frac{1}{2},\frac{3}{2},0,0),$ $(\frac{3}{2},\frac{1}{2},1,0),$ $(\frac{1}{2},\frac{3}{2},1,0),$ $(\frac{1}{2},\frac{1}{2},2,0),$ $(\frac{3}{2},\frac{3}{2},2,0).$ \\

\noindent $XZ$: $(\frac{1}{2},0,\frac{1}{2},0),$ $(\frac{1}{2},0,\frac{3}{2},0),$ $(\frac{3}{2},0,\frac{1}{2},0),$ $(\frac{1}{2},1,\frac{3}{2},0),$ $(\frac{3}{2},1,\frac{3}{2},0),$ $(\frac{1}{2},2,\frac{1}{2},0),$ $(\frac{3}{2},2,\frac{1}{2},0),$ $(\frac{3}{2},2,\frac{3}{2},0).$\\

\noindent $YZ$: $(0,\frac{1}{2},\frac{1}{2},0),$ $(0,\frac{1}{2},\frac{3}{2},0),$ $(0,\frac{3}{2},\frac{1}{2},0),$ $(1,\frac{1}{2},\frac{3}{2},1),$ $(1,\frac{3}{2},\frac{3}{2},1),$ $(2,\frac{1}{2},\frac{1}{2},0),$ $(2,\frac{3}{2},\frac{1}{2},0),$ $(2,\frac{3}{2},\frac{3}{2},0).$\\

\noindent $YW$: $(1,\frac{1}{2},1,\frac{1}{2}),$ $(1,\frac{1}{2},2,\frac{1}{2}),$ $(1,\frac{3}{2},1,\frac{1}{2}),$ $(1,\frac{3}{2},2,\frac{1}{2}).$\\

\noindent $ZW$: $(1,0,\frac{3}{2},\frac{1}{2}),$ $(1,2,\frac{3}{2},\frac{1}{2}).$

\noi These are explictely the 30 squares in $\R^4$ corresponding to the disk at the top right of figure 5 (after identifying diametrically opposite points).

\begin{rem} 
If $S$ is a gridded surface contained in the scaffolding $\mathcal{C}^2$ of the canonical cubulation $\cal C$  of $\mathbb{R}^{3}$, then 
$S$ is a gridded surface contained in the scaffolding $\mathcal{C}^2$ of the canonical cubulation $\cal C$  of $\mathbb{R}^{4}$, since $\mathbb{R}^{3}$
is canonically isomorphic to the hyperplane ${\cal{P}}=\mathbb{R}^3\times\{0\}\subset\mathbb{R}^4$ and  ${\cal{P}}$ has a canonical cubulation ${\cal{C}}_{\cal{P}}$ given by the restriction of the cubulation $\cal C$ of $\mathbb{R}^{4}$ to it; {\it i.e.}, ${\cal{C}}_{\cal{P}}$ is the decomposition into cubes which are the images of the unit cube
$I^{3}=\{(x_{1},x_2,x_{3},0)\,|\,0\leq x_{i}\leq 1\}$ by translations by vectors with integer coefficients whose last coordinate is zero. 
\end{rem}

\begin{figure}[h]  
\begin{center}
\includegraphics[height=3.5cm]{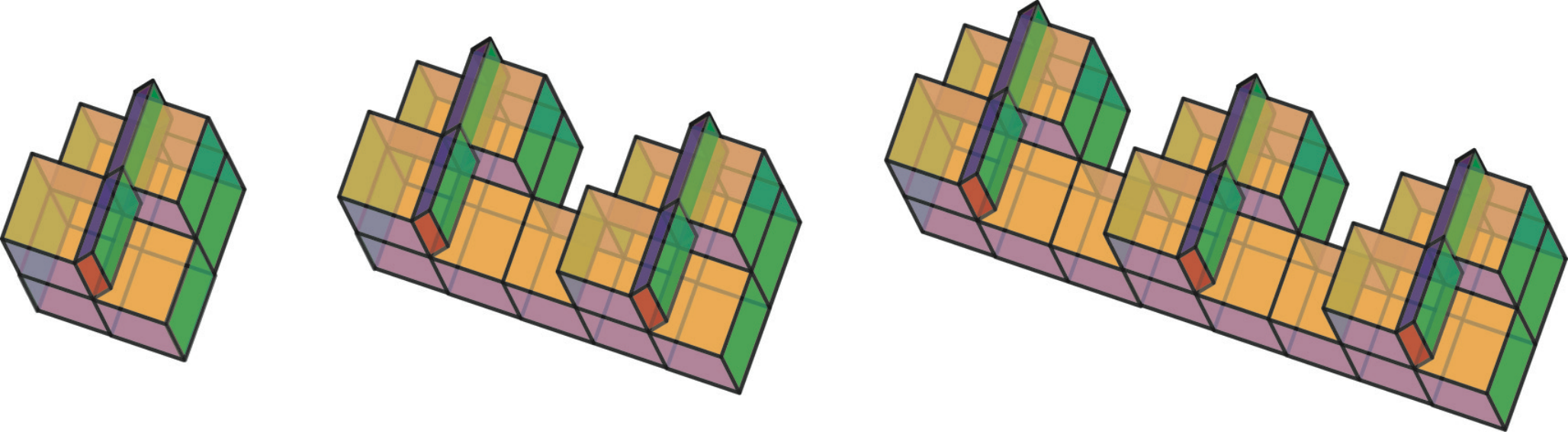}
\end{center}
\caption{\sl Closed gridded non-orientable surfaces in $\mathbb{R}^4$ (the first figure is the projection of the crosscap in 
$\R^3$ following the description of figure 5).} 
\label{F5}
\end{figure}

\begin{figure}[h]  
\begin{center}
\includegraphics[height=3.5cm]{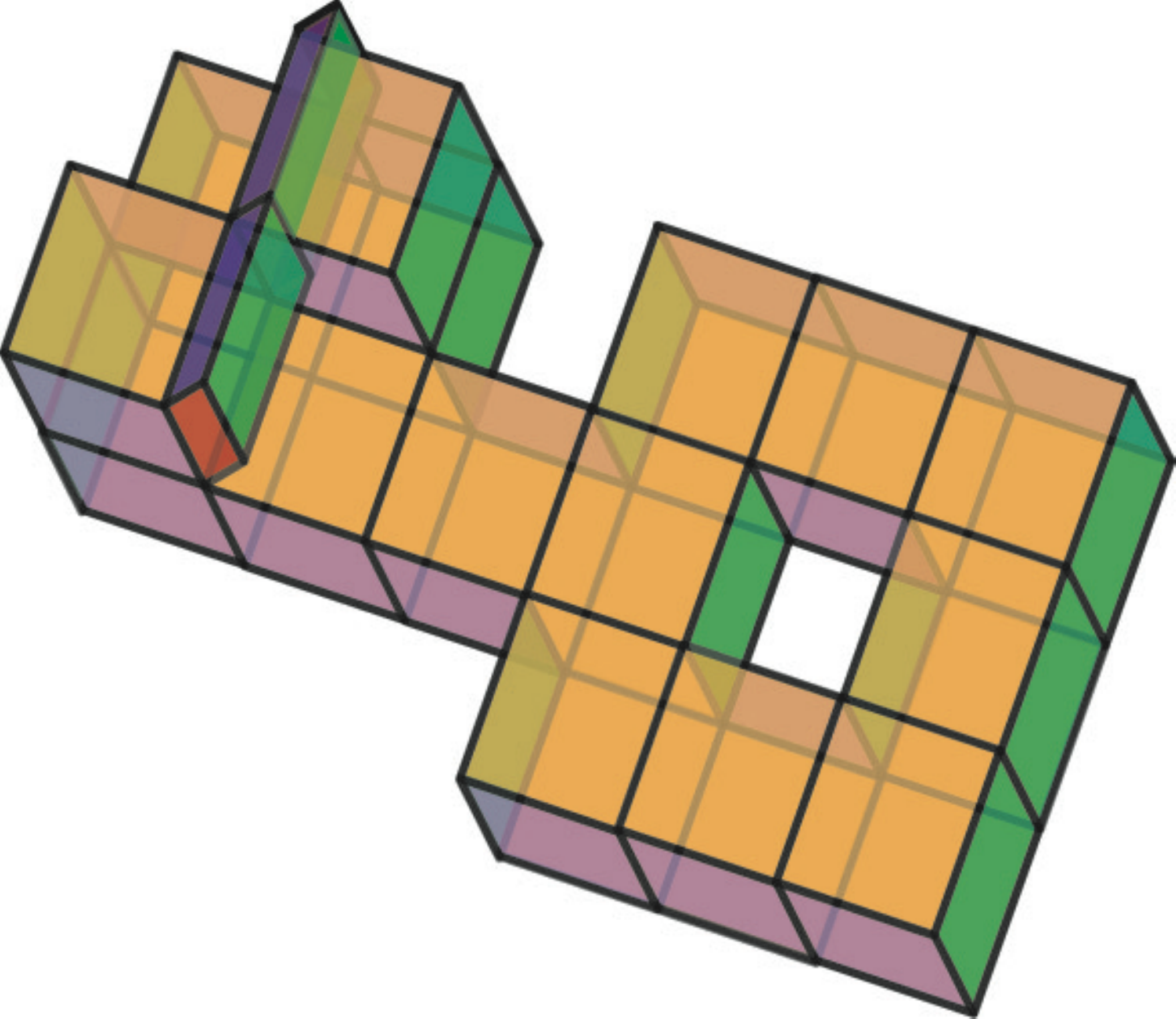}
\end{center}
\caption{\sl  The gridded connected sum of gridded surfaces is gridded.} 
\label{CS}
\end{figure}

\noindent Let $S_1$ and $S_2$ be two gridded surfaces. In a natural way, we define the \emph{gridded connected sum} $\#_g$ of $S_1$ and $S_2$, denoted by $S_1\#_g S_2$,  
as follows: We choose embeddings $i_j:\mathbb{D}^2\rightarrow S_j$ $(j=1,2)$, such that 
$i_j(\mathbb{D}^2)$ is a unit square $F_i$ into $S_i$, $i=1,2$. We can assume, up to  applying rigid movements that $F_1$ and $F_2$ are faces of some 3-cube $Q$ which 
its interior does not intersect neither $S_1$ or $S_2$. Thus we obtain
$S_1\#_g S_2$ from the disjoint sum $(S_1\setminus\Int (Q_1))\sqcup (S_2\setminus\Int (Q_2))$ joining $F_1$ with $F_2$ via the four remaining faces of $Q$  (see Figure \ref{CS}).
Observe that $S_1\#_g S_2$ is homeomorphic to the usual connected sum $S_1\# S_2$.\\

\begin{lem}
If $S_1$ and $S_2$ are surfaces homeomorphic to gridded surfaces $C_1$ and $C_2$, respectively; then the connected sum $S_1\# S_2$ is homeomorphic
to the gridded connected sum $C_1\#_g C_2$.
\end{lem}
\noindent \emph{Proof}. Consider the gridded surfaces $C_1$ and $C_2$. Then  $C_1\#_g C_2$ is homeomorphic to $C_1\# C_2$. Therefore 
$C_1\#_g C_2$ is homeomorphic to $S_1\# S_2$. $\square$\\

\noindent Summarizing, from the above Lemmas and using the classification theorem for closed surfaces, we have the following.
\begin{theo}
Any closed surface $S$ is homeomorphic to a gridded surface $C$ such that if $S$ is orientable then $C$ is contained in 
the scaffolding $\mathcal{C}^2$ of the canonical cubulation $\cal C$  of $\mathbb{R}^{3}$, and if $S$ is non orientable then $C$ is contained in the scaffolding $\mathcal{C}^2$ of the canonical cubulation $\cal C$  of $\mathbb{R}^{4}$.
\end{theo}

\subsection{Gridded surfaces in $\{4,3,4\}$ and $\{4,3,3,4\}$ in $\mathbb{R}^3$ and $\mathbb{R}^4$}

\noi We will show that all connected surfaces (orientable or not) can be gridded in $\mathbb{R}^4$ and moreover all orientable connected surfaces can be gridded 
in $\mathbb{R}^3$.\\ 

\noindent This follows from the fact that the 3-regular infinite tree graph can be embedded in the 1-skeleton of the 4-regular tessellation of squares in the Euclidean plane
(canonical cubulation of  $\mathbb{R}^2$). Then we thicken this tree graph to get a gridded tree of life in 
$\mathbb{R}^3$ or $\mathbb{R}^4$. In order to obtain a connected surface from our gridded tree of life,  we 
glue a set of tori and/or projective planes to the pruned tree of life corresponding to the given surface.

\begin{lem} 
The 3-regular infinite tree graph (\emph{i.e.,} a connected infinite regular tree of degree 3)  can be embedded in the 1-skeleton of the 4-regular tessellation of squares in the
 Euclidean plane (canonical cubulation of $\mathbb{R}^2$).
\end{lem}

\begin{figure}[h]  
\begin{center}
\includegraphics[height=3.5cm]{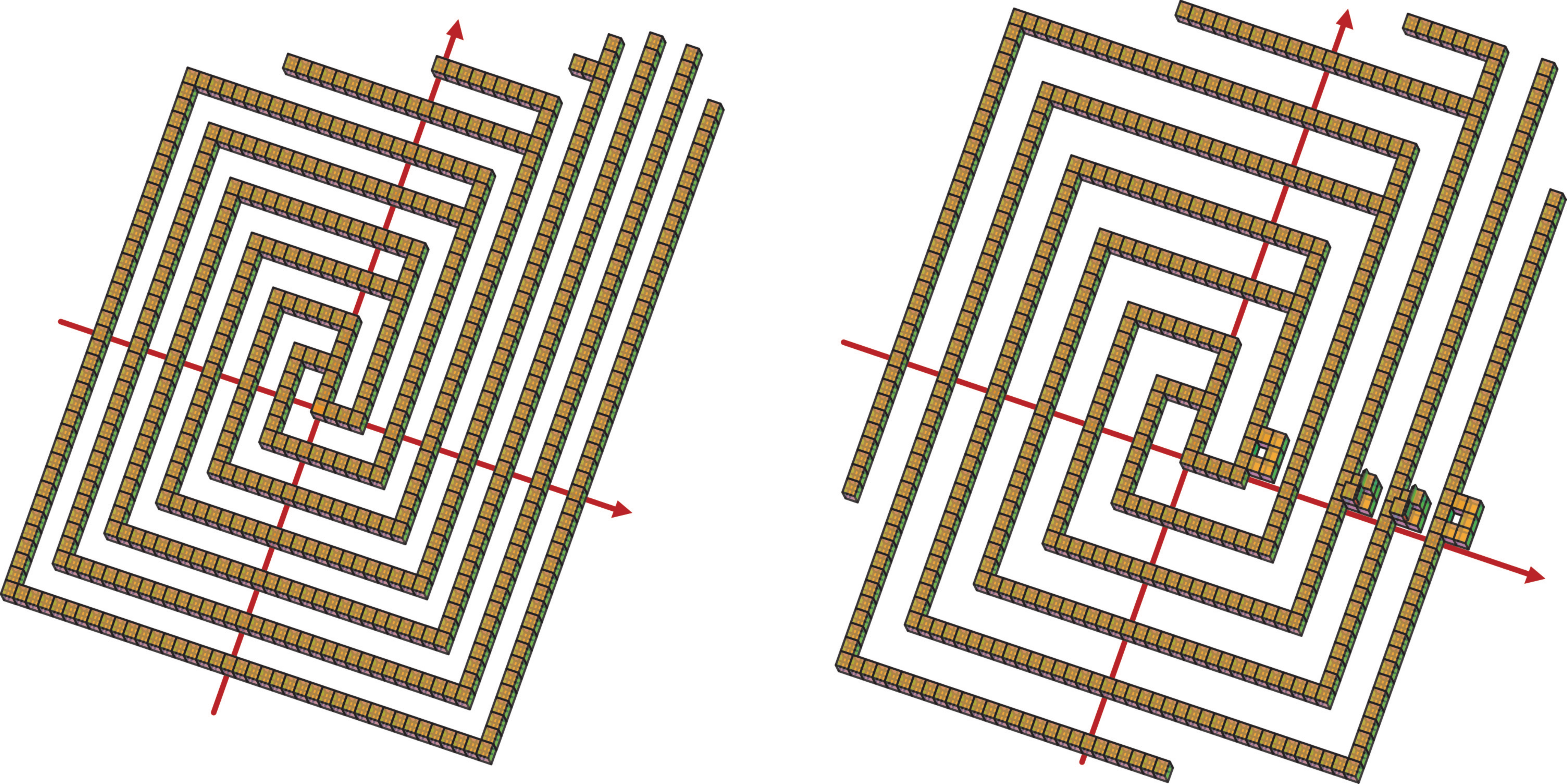}
\end{center}
\caption{\sl Left: The tree of life gridded in $\mathbb{R}^{3}$. Right: A general noncompact surface gridded in $\mathbb{R}^{4}$.} 
\label{TG}
\end{figure}

 \proof We will construct a 3-regular infinite tree graph embedded in the 1-skeleton of the 4-regular tessellation of squares in the Euclidean plane 
 $\mathbb{R}^2$ as a spiral. The trivalent vertices are the elements of the set $\{(n,2n) | n \in \mathbb{N}\}$ (see left part of Figure \ref{TG}). More precisely, we start from the origin $O=(0,0)$, 
and we consider the following two paths: one from $O$ to $(0,1)$ and the other one from $O$ to $(1,0)$ 
 and then from $(1,0)$ to $(1,3)$. Next we set two paths $\gamma^n_1$ and $\gamma^n_2$
  at each vertex $(n,2n+1)$, $n\in\mathbb{N}$. Notice that each path in the 1-skeleton of our cubulation of $\mathbb{R}^2$ can be described in a 
  unique way by a sequence of adjacent vertices on the $\mathbb{Z}^2$-lattice which is described by arrows. Using the above, we have that 
$$
{\scriptstyle\gamma^n_1=[(n,2n+1)\to(-2n-1,2n+1)\to(-2n-1,-2n-1)\to(2n+2,-2n-1)\to(2n+2,4n+5)]}  
$$
and 
$$
{\scriptstyle\gamma^n_2=[(n,2n+1)\to(n,2n+2)\to(-2n-2,2n+2)\to(-2n-2,-2n-2)\to(2n+3,-2n-2)to(2n+3, 4n+7)].}
$$

\noi The infinite union of paths is our desired 3-regular infinite tree graph. $\square$

\begin{theo} 
Any connected surface $S$ is homeomorphic to a gridded surface $C$ such that if $S$ is orientable then $C$ is contained in 
the scaffolding $\mathcal{C}^2$ of the canonical cubulation $\cal C$  of $\mathbb{R}^{3}$, and if $S$ is non orientable then $C$ is contained in the scaffolding $\mathcal{C}^2$ of the canonical cubulation $\cal C$  of $\mathbb{R}^{4}$.
\end{theo}

 \proof Given a connected surface $S$, we will construct a homeomorphic gridded copy of it  by means of Richard's Theorem. First, we recall that $S$ can be obtained from
 the tree of life by removing a set of branches to get the corresponding pruned tree  of life $P$ and next, we interchanged a finite or infinite number of pair of pants of $P$ 
 by either the connected sum of a pair of pants with a torus 
or the connected sum of a pair of pants with a projective plane to get $S$ (see Theorem \ref{clas}).\\
 
\noi Let $T$ be the infinite regular tree graph of degree 3. Then by the previous Lemma, $T$ can be embedded in the 1-skeleton of the 4-regular tessellation by squares in the 
Euclidean plane $\mathbb{R}^2$. Now, we apply the homothetic transformation $x\mapsto5x$ in $\R^2$ to expand 
$T$ obtaining a new tree graph $\hat{T}$ whose edge size is now 5 units.  Notice that the plane $\R^2$ is embedded in a natural way into $\R^3$ as the set of points whose third coordinate is zero. Let $V$ be the union of all cubes that intersect $\hat{T}$, so its boundary $\mathfrak{T}$ is a surface homeomorphic to the tree of life
embedded in the 2-skeleton of the canonical cubulation of $\R^3$ (see Figure \ref{TG}). By Richard's Theorem, we obtain $S$ gluing  a set of tori and/or projective planes 
to the pruned tree of life corresponding to $S$.
Notice that the space of ends $E_{\mathfrak{T}}$ of the tree of life $\mathfrak{T}$ is homeomorphic to the Cantor set, $C$, and our surface $S$ is determined by a 
nested sequence of three closed subsets of the Cantor set  $(E_S, G_S, O_S)$ where
 $O_S\subset{G_S}\subset{E_S}$. \\

\noi Let $\phi: E_{\mathfrak{T}} \rightarrow C$ be a homeomorphism of the set of ends of our tree of life into the Cantor set, and consider a homeomorphism $\Psi:E_S\to E'\subset{E_{\mathfrak{T}}}$ of the sets of ends of $S$ into the set of ends of our tree of life. The ends that do not belong to $\Psi(E_S)$ are pruned and closed by squares to obtain a surface without boundary. Observe that the pair of pants of the tree of life are in correspondence with the vertices $(0,n), n\in \mathbb{N}$, hence we glue tori and projective planes into pair of paints around these vertices $(0,n)$ and according to the images of $(E_S, G_S, O_S)$ into $E_{\mathfrak{T}}$. First glue tori on the respective pair of paints of the orientable non planar ends $O_S$ an next glue projective planes on the respective pair of paints of the ends in $G_S - O_S$ (see right part of Figure \ref{TG}).

\noi If the surface $S$ is non orientable but $G_S - O_S=\empty$ the there are two cases: If $O_S\neq \empty$ only is necessary glue one or two projective planes if $S$ is non orientable non or even, respectively. If $O_S= \empty$ then its necessary glue a finite number of projective planes and/or tori.   

\noi Notice that they are disjoint since the diameter of such pieces is less than the distance of the pair of parts of our surface. Therefore, $S$ is gridded in $\mathbb{R}^4$ or $\mathbb{R}^3$ if $S$ is orientable. $\square$

\begin{figure}[h]  
\begin{center}
\includegraphics[height=4cm]{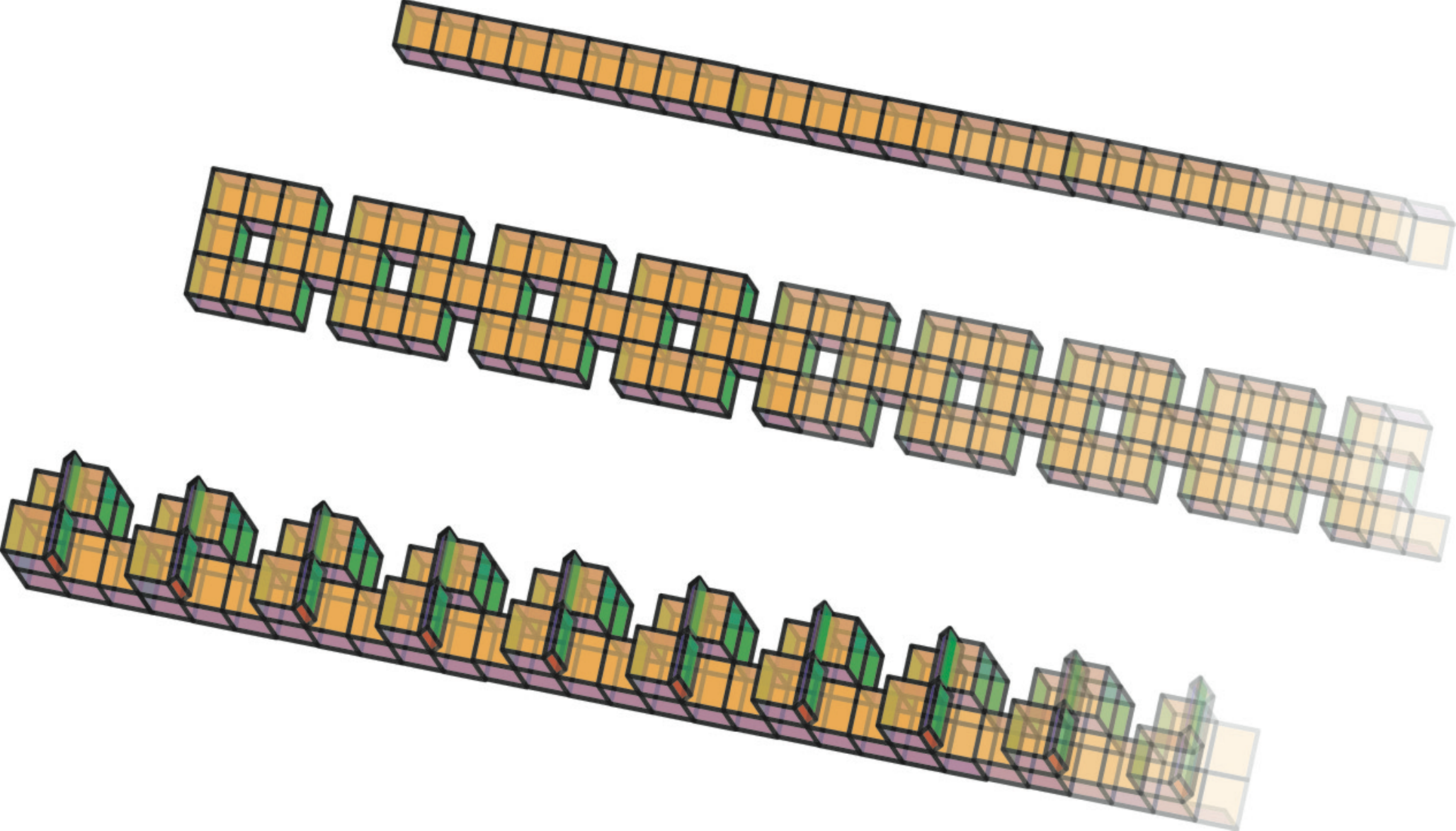}
\end{center}
\caption{\sl Types of ends of a surface with finite ends. A cylinder, an infinite ladder and an infinite chain of crosscaps.} 
\label{ends}
\end{figure}

\begin{rem} If the surface to considerer has a finite number of ends the we can consider a closed gridded surface and glue the ends which are of three types: the cilinder, the infinite ladder (an infinite chain of tori) and an infinite chain of projective planes. See Figure \ref{ends}.
\end{rem}

%\begin{figure}[h]  
%\begin{center}
%\includegraphics[height=3.5cm]{eucnoc.pdf}
%\end{center}
%\caption{\sl A gridded orientable non compact surfaces.} 
%\label{F5}
%\end{figure}

\section{Orientable surfaces in $\{4,3,5\}$ in $\mathbb{H}^3$}
 
All orientable surfaces can be constructed as gridded surfaces on the hyperbolic cubic honeycomb $\{4,3,5\}$ of the hyperbolic space $\mathbb{H}^3$. We proceed as in the Euclidean case where we proof that the orientable 
closed surfaces are gridded in $\mathbb{R}^3$ by means gridded the torus and the connected sum of two gridded surfaces.

\begin{figure}[h]  
\begin{center}
\includegraphics[height=2cm]{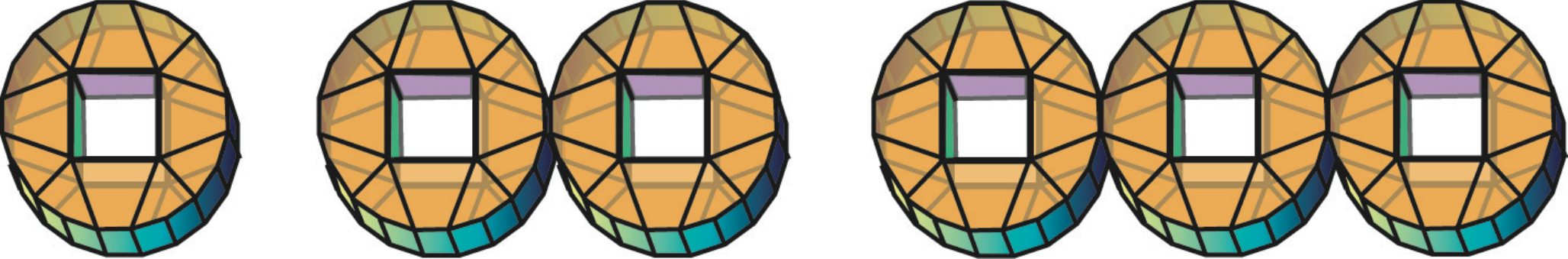}
\end{center}
\caption{\sl Torus and closed orientable surfaces in $\{4,3,5\}$ in $\mathbb{H}^3$.} 
\label{hyptorus}
\end{figure}

\begin{lem}
The torus is a gridded surface in $\{4,3,5\}$ in $\mathbb{H}^3$. Moreover, all closed orientable surfaces are gridded surfaces in $\{4,3,5\}$ in $\mathbb{H}^3$. 
\end{lem} 
 \proof
\noindent Let the torus be the gridded surface obtained as the boundary of twelve consecutive cubes in $\{4,3,5\}$ whose union  looks like as O. There are a \textit{central removed cube} and there are 4 cubes around each of its 
four hyperparallel edges (see Figure \ref{hyptorus}). This gridded torus in hyperbolic space is the one with the minimum number of squares in the scaffolding of $\{4,3,5\}$. It has 44 squares.

\noi There is a completely analogous hyperbolic concept of connected sum for gridded surfaces as in the previous Euclidean section. Let $S_1$ and $S_2$  gridded surfaces in $\mathbb{H}^3$ and 
$D_1\subset S_1$ and $D_2\subset S_2$ two squares such that each one of the corresponding support hyperbolic plane $\hat{D_i}$ $i=1,2$, divides $\mathbb{H}^3$ in two half--spaces in such a way that  $S_i$ is 
contained in only one half--space. Then we can construct the connected sum $S_1 \# S_2 $ as a gridded surface. A closed orientable surface can be gridded in this hyperbolic context as a connected sum of gridded tori. 
\endproof

\noi A \emph{hyperbolic pair of pants} is a closed pair of pants with a hyperbolic metric such that each of its three boundary circles are geodesics. The isometry class of such pair of pants is determined by the triple of lengths $(l_1,l_2,l_3)$ of the boundaries.

\begin{lem}
The pair of pants is a gridded surface in $\{4,3,5\}$ in $\mathbb{H}^3$. 
\end{lem} 

 \begin{figure}[h]  
\begin{center}
\includegraphics[height=3.5cm]{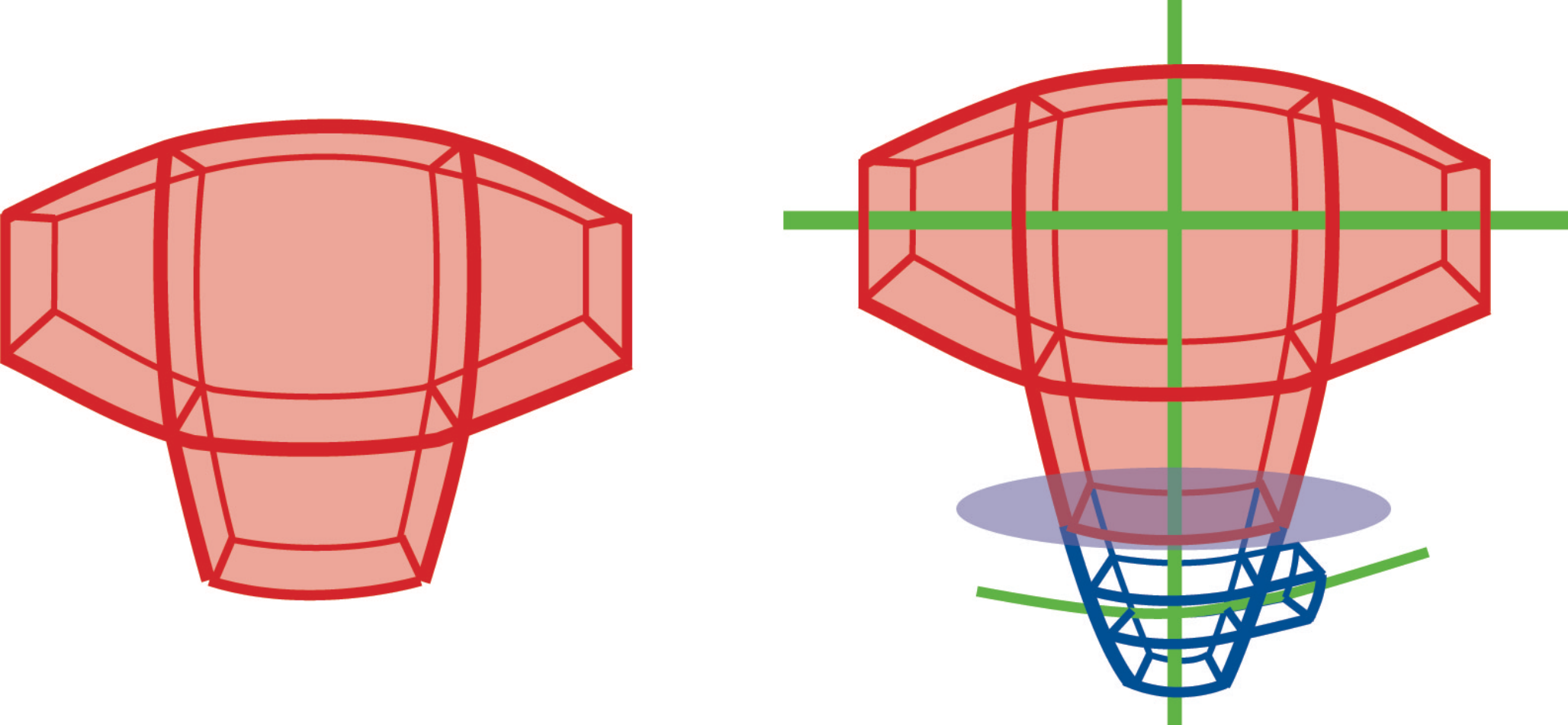}
\end{center}
\caption{\sl A gridded pair of pants in $\{4,3,5\}$ in $\mathbb{H}^3$.} 
\label{pairofpants}
\end{figure}

 \proof
\noindent Let the pair of pants be the gridded surface obtained as the boundary of four cubes in $\{4,3,5\}$ whose union looks like as T minus three squares. There are a \textit{central cube} $C$ and three neighborhood cubes $C_i$
of it such that these four cubes do not share a vertex; \emph{i.e.} $C\cap \cup_{i=1}^3 C_i=\emptyset$. Then $C\#C_1\#C_2\#C_3$ is our pair of pants (see Figure \ref{pairofpants}).
\endproof

 \begin{lem}
The tree of life is a gridded surface in $\{4,3,5\}$ in $\mathbb{H}^3$. 
\end{lem}

 \proof The proof is constructive by means pair of pants pasted along their boundaries as an infinite tree.\\
 
\noindent There are two distinguished geodesics in our model of the pair of pants. Notice that a pair of pants has a rotational symmetry of order 2. The \textit{axis of symmetry} is a geodesic which pass through the barycenters of the central cube 
and the second neighborhood cube \emph{i.e.} the vertical bar in the T. The \textit{second axis} is the geodesic perpendicular to the axis of symmetry which passes through the barycenters of the central cube and the first 
and third neighborhood cubes \emph{i.e.} the horizontal bar in the T.   
The axis of symmetry is ultraparallel to the two squares which were removed from the first and third neighborhood cubes and the second axis is ultraparallel to the square which was removed from the second neighborhood cube. \\

\noindent Then we can construct inductively the tree of life. We start by a pair of pants $P$ and glue it both a square in the boundary of the second cube and two pairs of pants $P_i$ 
at each boundary of the first and third cubes of $P$ with the second neighborhood cubes of the corresponding $P_i$; in such away that all barycenters of the cubes are lie in one hyperbolic plane (see Figure \ref{pairofpants}). 
The second axis of the first pair of pants is the axis of symmetry of the two pair of pants which have been glued to  its first and third neighborhood cubes. 
The second axis of the new two cubes is ultraparallel to the axis of symmetry of the original cube and these geodesics are ultraparallel to the hyperbolic plane defined by the square of gluing of each two pair of pants. 
This plane divides the tree of life in two connected components. \\

 \begin{figure}[h]  
\begin{center}
\includegraphics[height=3.5cm]{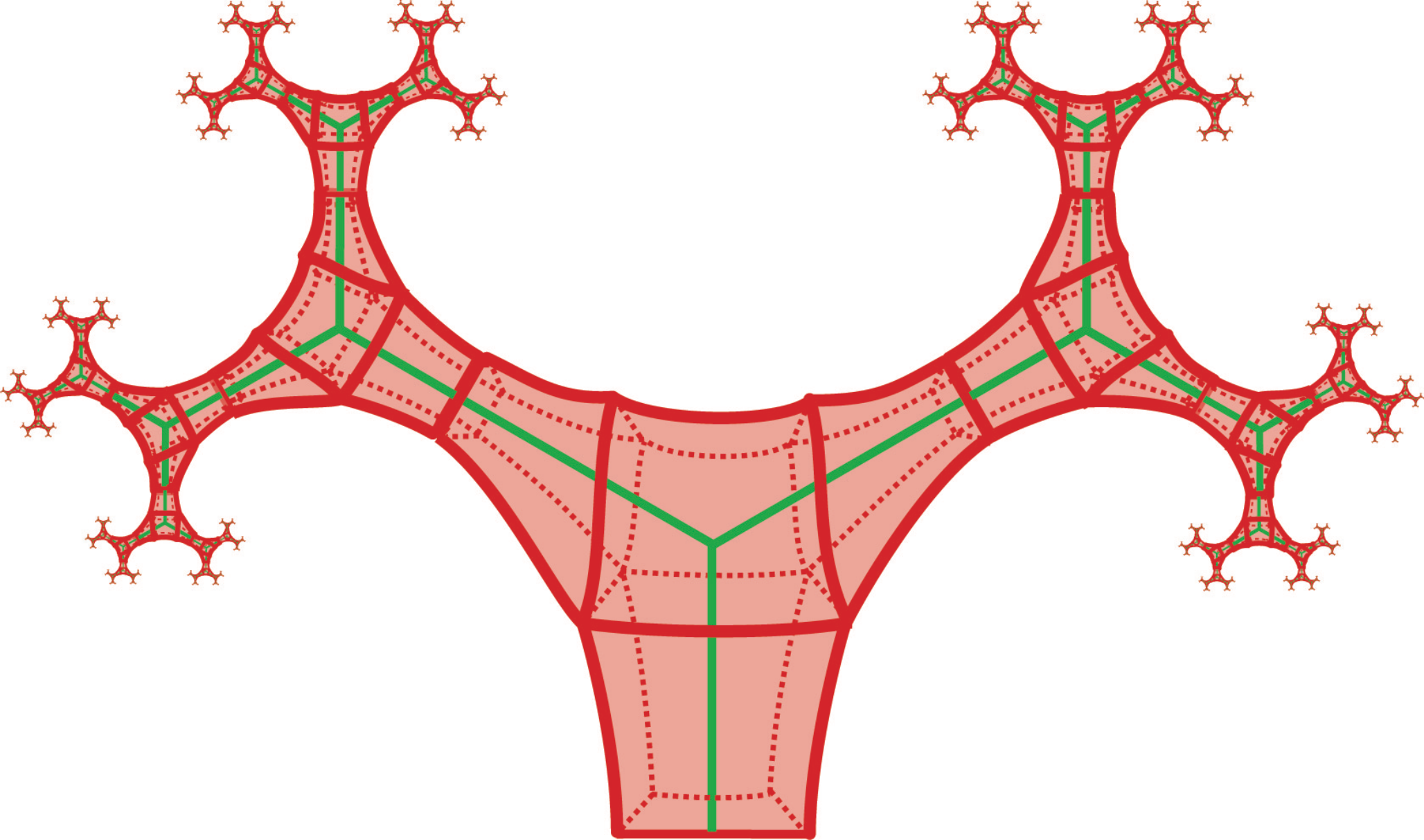}
\end{center}
\caption{\sl Gridded Tree of life in $\{4,3,5\}$ in $\mathbb{H}^3$.} 
\label{treeoflife}
\end{figure}

\noindent We glue new pairs of pants by the second neighborhood cube in the boundaries of the surface such that all barycenters of the cubes are in the hyperbolic plane. 
The step of induction is that the second axis of the pair of pants in the surface is the axis of symmetry of the two pair of pants which have been glued to their first and third neighborhood cubes. 
The second axis of the new two cubes is ultraparallel to the axis of symmetry of the original cube and these geodesics are ultraparallel to the hyperbolic plane defined by the square of gluing of each two pair of pants. 
This plane divides the tree of life in two connected components. By induction, we have constructed the tree of life.
\endproof 
\noindent We are ready to prove the following theorem.

\begin{theo}\label{os}
Any connected orientable surface is homeomorphic to a gridded surface in $\{4,3,5\}$ in $\mathbb{H}^3$. 
\end{theo}

\noindent \emph{Proof}. Any connected orientable surface can be constructed from the pruned tree of life replacing some pair of pants by ``handles''. 
We can consider the gridded torus minus three nonconsecutive equatorial squares. Notice that when we exchange a pair of pants by these gridded tori minus three 
nonconsecutive equatorial squares, the property of the gridded connected sum is preserved. 
The hyperbolic planes which pass through the boundaries of the modified pair of pants are ultraparallels, then the construction of a tree of life with handles is analogous to the planar tree of life.

\section{Surfaces in $\{4,3,3,5\}$ in $\mathbb{H}^4$}

\noindent One great difference between the Euclidean and the hyperbolic gridded cases is that the gridded Euclidean spaces are nested and the gridded hyperbolic spaces are not (as gridded spaces). For the Euclidean case we proof that the orientable closed surfaces are 
gridded in $\mathbb{R}^3$ by means gridded the torus and the connected sum. We proved that the projective plane is gridded in $\mathbb{R}^4$ and all closed surfaces are gridded in $\mathbb{R}^4$.

\noindent The gridded hyperbolic 3-space $\{4,3,5\}$ is not contained in the gridded hyperbolic 4-space $\{4,3,3,5\}$. 
However, it is not a great problem to grid in $\mathbb{H}^4$ all the gridded surfaces in $\mathbb{H}^3$. We need to grid the torus, the pair of pants and the connected sum of gridded surfaces in order to obtain 
all orientable surfaces. For nonorientable surfaces we need only to prove that the projective plane is gridded  in $\mathbb{H}^4$ and applying the Richards Theorem then we will obtain all surfaces gridded in $\mathbb{H}^4$.

\begin{figure}[h]  
\begin{center}
\includegraphics[height=2.5cm]{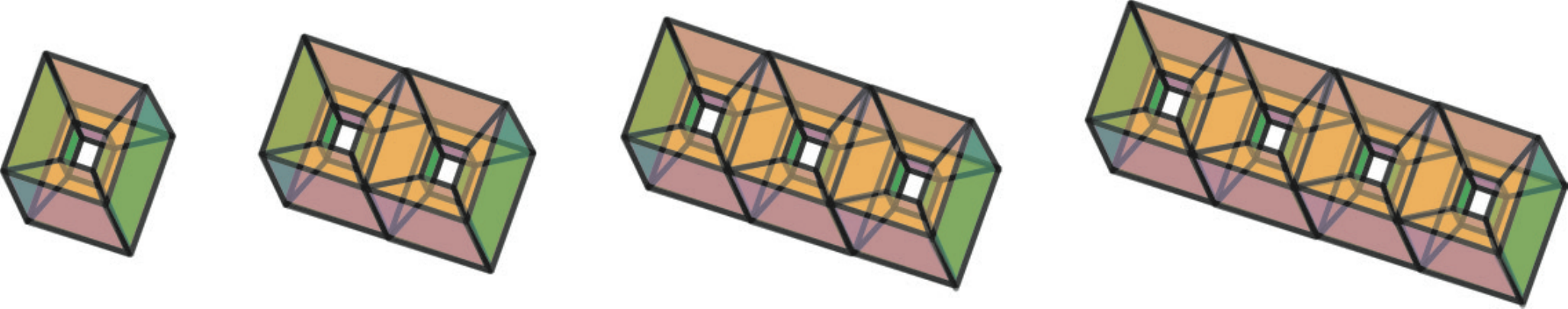}
\end{center}
\caption{\sl Gridded torus and closed orientable surfaces in $\{4,3,3,5\}$ in $\mathbb{H}^4$.} 
\label{torush4}
\end{figure}

\begin{lem}
The torus is a gridded surface in $\{4,3,3,5\}$ in $\mathbb{H}^4$. Moreover, all closed orientable surfaces are gridded surfaces in $\{4,3,3,5\}$ in $\mathbb{H}^4$. 
\end{lem} 
 \proof
\noindent There are 8 cubes in the hypercube forming  two linked tori in $\mathbb{S}^3$ (see Figure \ref{torush4}). We take one of these torus.
Let the torus be the gridded surface obtained as the boundary of four consecutive cubes in a hypercube $\{4,3,3\}$.

\noi As above, we have a 4-dimensional analogous hyperbolic concept of connected sum for gridded surfaces. Let $S_1$ and $S_2$ be two hyperbolic gridded surfaces $\mathbb{H}^4$ and 
$D_i\subset S_i$, $i=1,2$ be squares such that the corresponding support hyperbolic plane $\hat{D_i}$ lies in a 3-dimensional geodesic hyperbolic subspace which divides $\mathbb{H}^4$ in two half--spaces 
such that $S_i$ is contained in only one half--space. Then we can construct the connected sum $S_1 \# S_2 $ as a gridded surface. A closed orientable surface can be gridded in this hyperbolic context 
as a connected sum of gridded tori. 
\endproof

\begin{lem}
The pair of pants is a gridded surface in $\{4,3,3,5\}$ in $\mathbb{H}^4$. 
\end{lem}

\begin{figure}[h]  
\begin{center}
\includegraphics[height=2cm]{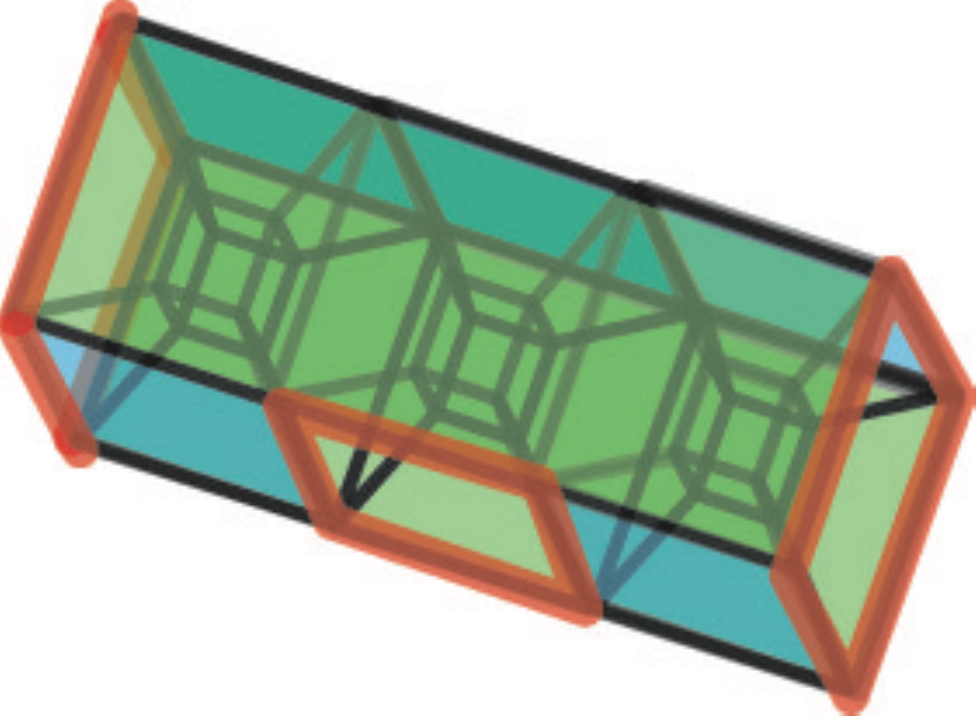}
\end{center}
\caption{\sl The gridded pair of pants in $\{4,3,3,5\}$ in $\mathbb{H}^4$.} 
\label{pair4}
\end{figure}

 \proof
\noindent Consider three 3-faces $F_1$, $F_2$ and $F_3$ of three consecutive hypercubes $C_1$, $C_2$ and $C_3$ such that $F_i\subset C_i$ ($i=1,2,3$) and  $F_1\cap F_2$ and $F_2\cap F_3$ are two disjoint squares. 
The barycenters of $C_1$, $C_2$ and $C_3$ are collinear. 
Consider the connected sum $F_1\#F_2\#F_3$. Notice that $F_1\#F_2\#F_3$ is homeomorphic to a $\mathbb{S}^2$. Then the pair of pants is obtained from $F_1\#F_2\#F_3$ by removing one boundary square $S_i$ from $F_i$,
in such a way that $S_1$ is parallel to $F_1\cap F_2$ and $S_3$ is parallel to $F_2\cap F_3$ (see Figure \ref{pair4}).
\endproof

 \begin{lem}
The tree of life is a gridded surface in $\{4,3,3,5\}$ in $\mathbb{H}^4$. 
\end{lem}

 \proof The proof is constructive using pair of pants pasted along their boundaries as an infinite tree.\\
 
\noindent There are two important geodesics in our model of a pair of pants. A pair of pants has a rotational symmetry of order 2. 
The \textit{axis of symmetry} is a geodesic which passes through the barycenter of $F_2$ and the barycenter of the square $S_2$.  
The \textit{second axis} is the geodesic perpendicular to the axis of symmetry which passes through the barycenters of $S_1$ and $S_3$.
The axis of symmetry is ultraparallel to the two squares which were removed from $F_1$ and $F_3$ and the second axis is 
ultraparallel to the square which was removed from $F_2$. \\

\noindent We can construct inductively the tree of life in analogous way as $\mathbb{H}^3$. \endproof

\begin{lem}
The projective plane is a gridded surface in $\{4,3,3,5\}$ in $\mathbb{H}^4$. 
\end{lem}

\begin{figure}[h]  
\begin{center}
\includegraphics[height=3.5cm]{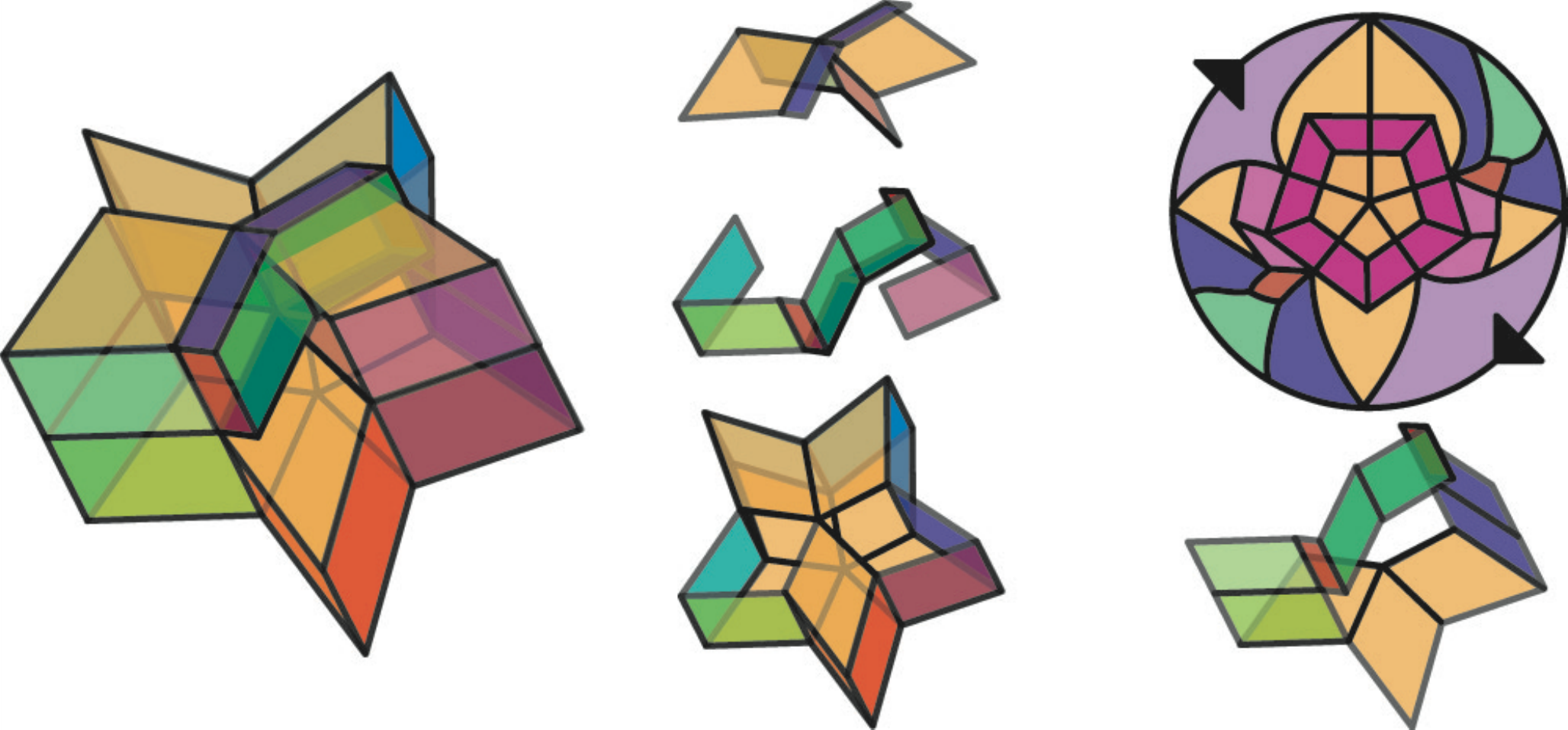}
\end{center}
\caption{\sl Gridded projective plane in $\{4,3,3,5\}$ in $\mathbb{H}^4$.} 
\label{PPh}

\end{figure}

\noindent \emph{Proof}. We construct a gridded version of the crosscap in $\{4,3,3,5\}$ in $\mathbb{H}^4$. 
See the Figures \ref{PP} and \ref{PPh}. In the left we show the projection of the crosscap. In the middle we divide the crosscap in three parts: 
at the bottom, there is the base which is a pentagonal cubic box minus two squares at the top. In the middle,
there is a band and at the top there is a disk which is a neighborhood of one vertex. Only the base is different in the two Figures  \ref{PP} and \ref{PPh}. \\ 

\noindent In the right part of the Figure \ref{PPh} there is a description of the combinatorial square complex of the crosscap as a  disk with 34 squares after identifying
points in the circle boundary by the antipodal map. At the bottom, there is a M\"obius band contained in this crosscap.
This is similar to the Euclidean case except one uses 4 more squares in the crosscap (see figure 5).

\noindent We are ready to prove the following theorem.

\begin{theo}
Any connected surface is homeomorphic to a gridded surface in $\{4,3,3,5\}$ in $\mathbb{H}^4$. 
\end{theo}

\begin{figure}[h]  
\begin{center}
\includegraphics[height=3.5cm]{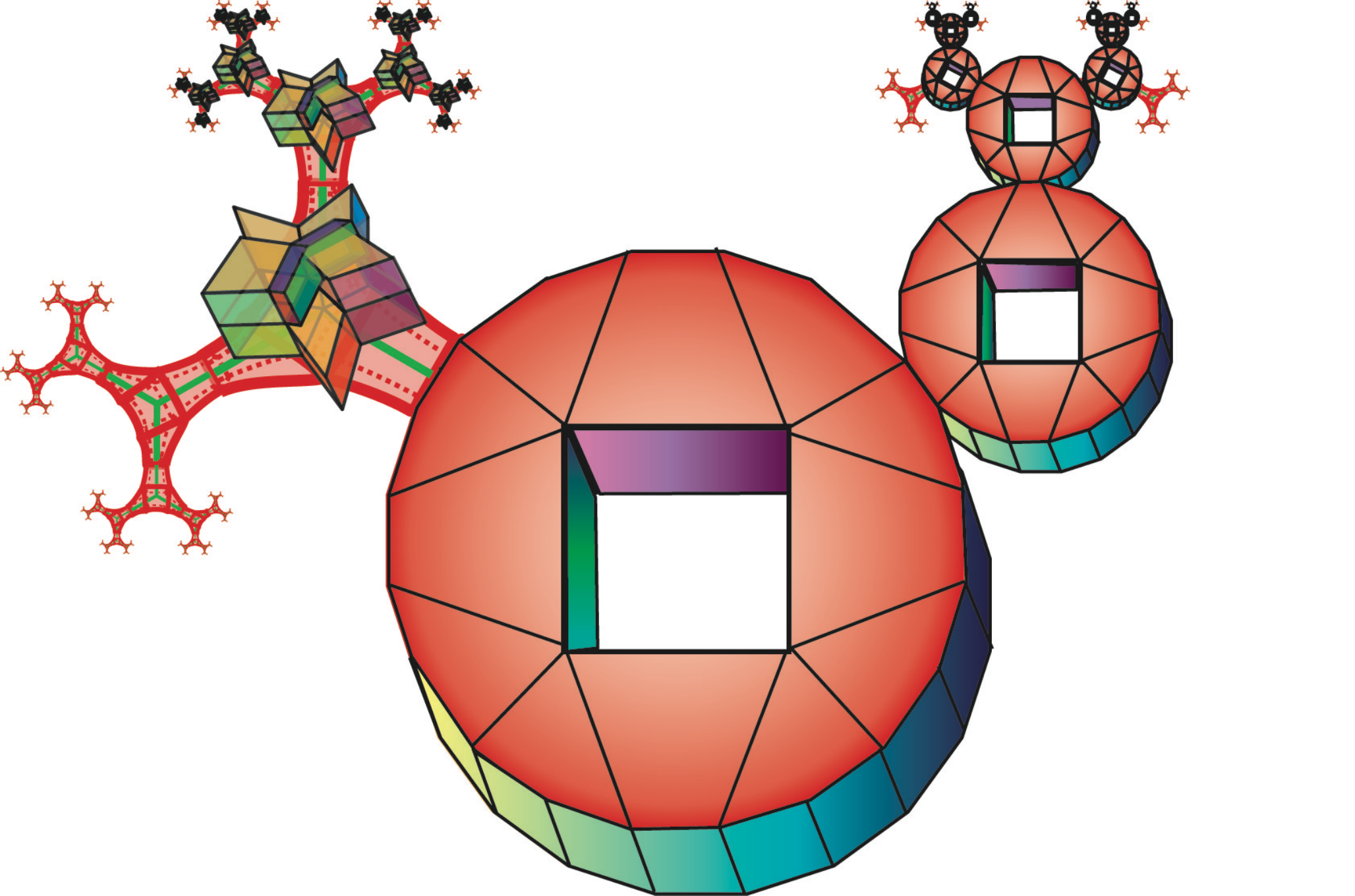}
\end{center}
\caption{\sl Tree of life with handles and projective planes. The non-planar non orientable ends go to the top left and the orientable non-planar ends go to the top right.} 
\label{ph}
\end{figure}

\noindent \emph{Proof}. Any connected surface can be constructed from the pruned tree of life 
modifying some pair of pants by means put ``handles'' and ``projective planes''. 
Consider the connected sum of a hyperbolic gridded torus in $\{4,3,3,5\}$ and the pair of pants. In fact, from a gridded torus we remove an open square and we  
paste the boundary of this square onto the boundary of a removed square in the central cube $F_2$ of the pair of pants (see Figure \ref{ph}). \\

\noi As above, we consider the connected sum of a hyperbolic gridded projective plane in $\{4,3,3,5\}$ and the pair of pants.
Remove from the gridded projective plane an open square in its base and paste its boundary with the boundary of one square in the central 
cube $F_2$ of the pair of pants. \\

\noi Notice that, if we exchange a pair of pants of a pruned tree of life  by these kind of new pair of pants the property of the gridded connected sum is preserved. 
The hyperbolic spaces which pass by the boundaries of these new pair of pants are ultraparallels and divide $\mathbb{H}^4$ in two half--spaces 
where the pair of pants is contained in one component. Then the construction of a tree of life with handles and projective planes is analogous to the construction
of an orientable noncompact surface in $\mathbb{H}^3$ (see Theorem \ref{os}).

\noindent J. P. D\'iaz. {\tt Instituto de Matem\'aticas, Unidad Cuernavaca}. Universidad Nacional Au\-t\'o\-no\-ma de M\'exico.
Av. Universidad s/n, Col. Lomas de Chamilpa. Cuernavaca, Morelos, M\'exico, 62209.

\noindent {\it E-mail address:} juanpablo@matcuer.unam.mx
\vskip .3cm

\noindent G. Hinojosa. {\tt Centro de Investigaci\'on en Ciencias}. Instituto de Investigaci\'on en Ciencias B\'asicas y Aplicadas. Universidad Aut\'onoma del Estado de Morelos. Av. Universidad 1001, Col. Chamilpa.
Cuernavaca, Morelos, M\'exico, 62209. 

\noindent {\it E-mail address:} gabriela@uaem.mx 

\vskip .3cm
\noindent A. Verjovsky. {\tt Instituto de Matem\'aticas, Unidad Cuernavaca}. Universidad Nacional Au\-t\'o\-no\-ma de M\'exico.
Av. Universidad s/n, Col. Lomas de Chamilpa. Cuernavaca, Morelos, M\'exico, 62209.

\noindent {\it E-mail address:} alberto@matcuer.unam.mx

\end{document}